\documentclass[pre,superscriptaddress,twocolumn]{revtex4-1}
\usepackage{amsfonts,amsmath,amssymb,mathtools,graphicx}
\usepackage[colorlinks=true,allcolors=blue]{hyperref}
\usepackage[svgnames,dvipsnames,x11names]{xcolor}
\usepackage{tikz}

\newcommand{\HH}{\mathbb{H}}

\newcommand{\ea}[1]{\langle #1 \rangle}

\newcommand{\e}{\varepsilon}

\newcommand{\Ric}{\mathrm{Ric}}
\newcommand{\vol}{\mathrm{vol}}
\renewcommand{\epsilon}{\varepsilon}
\DeclareMathOperator{\arctanh}{arctanh}

\begin{document}

\title{Ollivier-Ricci curvature convergence in random geometric graphs}

\author{Pim van der Hoorn}
\affiliation{Department of Mathematics and Computer Science, Eindhoven University of Technology, 5600 MB, Eindhoven, Netherlands}
\affiliation{Department of Physics, Northeastern University, Boston, Massachusetts 02115, USA}
\affiliation{Network Science Institute, Northeastern University, Boston, Massachusetts 02115, USA}

\author{William~J.~Cunningham}
\affiliation{Perimeter Institute for Theoretical Physics, Waterloo, Ontario N2L2Y5, Canada}
\affiliation{Department of Physics, Northeastern University, Boston, Massachusetts 02115, USA}
\affiliation{Network Science Institute, Northeastern University, Boston, Massachusetts 02115, USA}

\author{Gabor~Lippner}
\affiliation{Department of Mathematics, Northeastern University, Boston, Massachusetts 02115, USA}

\author{Carlo Trugenberger}
\affiliation{SwissScientific Technologies SA, CH-1204 Geneva, Switzerland}

\author{Dmitri Krioukov}
\affiliation{Department of Physics, Northeastern University, Boston, Massachusetts 02115, USA}
\affiliation{Network Science Institute, Northeastern University, Boston, Massachusetts 02115, USA}
\affiliation{Department of Mathematics, Northeastern University, Boston, Massachusetts 02115, USA}
\affiliation{Department of Electrical and Computer Engineering, Northeastern University, Boston, Massachusetts 02115, USA}

\begin{abstract}
Connections between continuous and discrete worlds tend to be elusive. One example is curvature.
Even though there exist numerous nonequivalent definitions of graph curvature, none is known to converge in any limit to any traditional definition of curvature of a Riemannian manifold. Here we show that Ollivier curvature of random geometric graphs in any Riemannian manifold converges in the continuum limit to Ricci curvature of the underlying manifold, but only if the definition of Ollivier graph curvature is properly generalized to apply to mesoscopic graph neighborhoods. This result establishes the first rigorous link between a definition of curvature applicable to networks and a traditional definition of curvature of smooth spaces.
\end{abstract}

\maketitle

\section{Introduction}

Curvature is one of the most basic geometric characteristics of space. The original definitions of curvature apply only to smooth Riemannian or Lorentzian manifolds, but there exist numerous extensions of curvature definitions applicable to graphs, simplicial complexes, and other discrete structures. For a variety of reasons, these extensions have recently seen a surge of interest in areas as diverse as network/data science and quantum gravity.

In network science, graph curvature is interesting in general since many real-world networks were found to possess different flavors of geometry~\cite{boguna2021network}. Measuring network curvature is then a way to learn what this geometry really is~\cite{lubold2020identifying}. On the application side, graph curvature was used to characterize congestion in telecommunication networks~\cite{narayan2011large,jonckheere2011euclidean}, detect cancer cells~\cite{sandhu2015graph}, predict COVID-19 spreading~\cite{souza2020using} and stock market fluctuations~\cite{akguller2018geodetic}, analyze robustness and other properties of the Internet, financial, and brain networks~\cite{ni2015ricci,sandhu2016ricci,tadic2018origin,tadic2019functional,farooq2019network}, and in a number of classic applications such as community inference~\cite{sia2019ollivier} and network embedding in machine learning~\cite{gu2019learning}. In a majority of these applications, graph curvature appears as a signature of a (pathological) event or property of interest. How reliable such a signature is for a task at hand depends on whether the used definition of graph curvature does what is intended---that it is indeed curvature, and not something else.

In general relativity, scalar Ricci curvature is a key player since it appears in the Einstein-Hilbert action whose least-action extremization leads to Einstein's equations~\cite{carroll2019spacetime}. Therefore, various extensions of Ricci curvature to discrete structures representing  ``quantum spacetimes'' have been considered. The first step in that direction was based on replacing continuum manifolds with simplicial complexes~\cite{ambjrn2012nonperturbative}. To that end, variations of the seminal definition of discrete curvature by Regge~\cite{regge1961general} were developed and utilized. More recently, a collection of more minimalistic structures was considered, such as causal sets~\cite{bombelli1987space,henson2009causal,surya2019causal,benincasa2010scalar,belenchia2016continuum,glaser2016hartlehawking,glaser2018finite,eichhorn2019spectral,eichhorn2019steps,cunningham2020dimensionally} and random graphs~\cite{trugenberger2017combinatorial,kelly2019assembly,kelly2021emergence,klitgaard2018introducing,klitgaard2018implementing,klitgaard2020quantum}. Yet another class of approaches aimed at explaining spacetime emergence from ``quantum bits'' appear in space-from-entanglement proposals such as ER=EPR or tensor networks~\cite{maldacena2013cool,raamsdonk2010building,cao2017space,bao2017sitter,swingle2018spacetime,orus2019tensor,cotler2019locality}. In all these examples, one can hope that a particular discrete structure represents Planck-scale gravitational physics, only if discrete curvature of this structure converges to Ricci curvature of classical spacetime in the continuum limit. If such convergence does not hold in a particular approach, then the approach does not agree with general relativity, so it cannot be considered realistic, regardless of how attractive it is in other respects.

Quite contrary to the continuum world of Riemannian and Lorentzian manifolds where there is no real ambiguity regarding the definition of curvature, there exist not one or two but quite many nonequivalent definitions of curvature applicable to graphs. We briefly review the most prominent definitions in the next section. 
Unfortunately, none of these definitions are known to converge in the continuum limit to any traditional curvature of any Riemannian manifold. Therefore, it is often unclear how one should interpret different measurements of different curvatures in different graphs and networks, raising question concerning reliability and generalizability of predictions based on such measurements. Similarly, the evasiveness of classical limit in some approaches to quantum gravity can be often linked to unclarities concerning whether discrete ``quantum'' curvature converges to a continuous classical one. These and other problems related to curvature convergence become particularly acute in view of growing evidence that different graph curvature definitions may disagree even about the sign of curvature in some paradigmatic graphs and networks~\cite{samal2018comparative,prokhorenkova2020global}.

Here we show that Ollivier curvature~\cite{ollivier2007ricci,ollivier2009ricci,ollivier2010survey} of random geometric graphs~\cite{gilbert1961random,dall2002random,penrose2003random} in any Riemannian manifold converges to Ricci curvature of the manifold in the continuum limit. Unfortunately, this convergence does not and cannot hold for the standard definition of Ollivier curvature of microscopic one-link neighborhoods in sparse unweighted graphs~\cite{jost2014olliviers}. Fortunately, we find a natural generalization of Ollivier curvature applicable to mesoscopic neighborhoods in unweighted and weighted graphs. This generalization is a not-previously-considered key ingredient that allows us to prove the curvature convergence. Our proofs hold only when the sizes of mesoscopic graph neighborhoods fall within certain windows that depend on graph density and on how graph edges are weighted. We also perform large-scale simulations that agree with our proofs and suggest that the convergence windows are actually much wider than the limits imposed by our proof techniques.

These are the first result of the sort, linking rigorously a definition of graph curvature to the traditional Ricci curvature of a Riemannian manifold. To the best of our knowledge, the closest, in spirit, previous results are Cheeger {\it et al.}'s proof~\cite{cheeger1984curvature} of Regge's seminal observation~\cite{regge1961general} that an angle-defect-based curvature of increasingly finer-grained simplicial triangulations of a manifold converges to its Ricci curvature, and the more recent demonstration~\cite{benincasa2010scalar,belenchia2016continuum} that the Benincasa-Dowker action of causal sets converges to its classical counterpart in the continuum limit.

To proceed, we first recall the Ollivier curvature definition in perspective of other curvature definitions applicable to graphs and networks, Sec.~\ref{sec:curvature}. We then describe our generalization of the Ollivier curvature definition to mesoscopic neighborhoods in graphs in general and in random geometric graphs in particular, Sec.~\ref{sec:mesoscopic}. We then state our main results in Sec.~\ref{sec:results}, and conclude with a discussion of their implications and caveats in Sec.~\ref{sec:discussion}. 

\section{Graph curvature definitions}~\label{sec:curvature}

There are many nonequivalent definitions of curvature applicable to graphs. It is impossible to list them all. Here we mention some notable ones due to Steiner~\cite{steiner1882uber}, Regge~\cite{regge1961general,cheeger1984curvature}, Bakry-\'Emery~\cite{bakry1985diffusions}, Gromov~\cite{gromov1987hyperbolic,jonckheere2008scaled}, Higuchi~\cite{higuchi2001combinatorial}, Eckmann-Moses~\cite{eckmann2002curvature}, Forman~\cite{forman2003bochners,sreejith2016forman,weber2018coarse}, Lin-Yau~\cite{lin2010ricci,lin2011ricci}, Knill~\cite{knill2010discrete,knill2012index}, Keller~\cite{keller2011curvature,keller2011cheeger}, and hybrids thereof~\cite{kempton2020large}. We first outline the taxonomy of these curvatures and basic ideas behind their definitions to put Ollivier curvature~\cite{ollivier2007ricci,ollivier2009ricci,ollivier2010survey} in perspective, and then recall the detailed definition of the latter for completeness. 

With a few exceptions, the discrete curvature definitions above fall into the following classes.

\subsection{Direct approaches}

The first class of definitions comes from relatively direct extensions of continuous curvature to nonsmooth objects such as simplices, simplicial complexes, or more general complexes. Since graphs are 1-dimensional simplicial complexes, any notion of curvature applicable to complexes applies to graphs as well. Casting cliques in graphs as higher-dimensional simplices turns graphs into higher-dimensional complexes, so that curvature definitions applicable to higher-dimensional complexes are also applicable to graphs.

One subclass of such approaches first defines curvature of polyhedra using the Gauss-Bonnet theorem that relates the curvature of a surface to its Euler's characteristic. This definition is then extended to general simplicial complexes. Curvatures defined by Steiner~\cite{steiner1882uber}, Regge~\cite{regge1961general,cheeger1984curvature}, Knill~\cite{knill2010discrete,knill2012index}, and Keller~\cite{keller2011curvature,keller2011cheeger} all fall into this category.

Forman curvature~\cite{forman2003bochners,sreejith2016forman,weber2018coarse} follows similar ideas, but is slightly different in that it uses the Bochner-Weitzenb\"{o}ck formula that relates the Laplacian of a Riemannian manifold to its discrete approximation, which is then used to derive a discrete version of curvature applicable to a very general class of complexes. The definition of Forman curvature for unweighted graphs is one of the simplest. If a graph is considered as a 1-complex, so that only nodes (0-simplices) and links (1-simplices) are considered, then the Forman curvature $F_1(i,j)$ of a link between nodes $i$ and $j$ of degrees $k_i$ and $k_j$ is
$
	F_1(i,j) = 4 - (k_i + k_j).
$
If the graph is considered as a 2-complex, so that triangles are also considered, and $\Delta_{i,j}$ is the number of triangles containing the link, then the definition becomes
$
	F_2(i,j) = F_1(i,j) + 3 \Delta_{i,j}.
$

\subsection{Indirect approaches}

Another class of ideas behind discrete curvature are less direct as they are based on more sophisticated curvature-related properties of smooth spaces. 

One such idea, pioneered by Bakry and Emery~\cite{bakry1985diffusions}, relies on the characterization of Ricci curvature via curvature-dimension inequalities.
The curvature-dimension inequality establishes a lower bound for the Ricci curvature of a Riemannian manifold based on the gradient of harmonic functions on it. The inequality is an immediate consequence of Bochner's identity that relates Ricci curvature to the gradient of harmonic functions.
It was an important insight by Bakry and Emery that the curvature-dimension inequality can be lifted from manifolds and considered as a substitute of the lower bound for Ricci curvature in spaces where direct discretizations of Ricci curvature are not possible. In addition to the Bakry-Emery definition~\cite{bakry1985diffusions}, the definitions of Lin-Yau~\cite{lin2010ricci,lin2011ricci} also fall into this category.

The key point behind this type of approaches to discrete curvature is to base its definitions on some fundamental nontrivial curvature-related properties that transcend well beyond the limits of smooth spaces.

\subsection{Ollivier curvature}

Ollivier curvature belongs to the second class of approaches outlined above. The fundamental aspect of Ricci curvature that Ollivier curvature embarks on to sail to the discrete world, is somewhat different from---although related to---the curvature-dimension inequalities. This aspect is of a less analytic and more geometric nature. It deals with how balls shrink or expand under parallel transport. Recall that if Ricci curvature of a space is negative, zero, or positive, then under parallel transport, balls in this space expand, stay the same, or shrink, respectively.

To capture this property, one needs to consider the Wasserstein a.k.a.\ transportation distance $W(\mu_1,\mu_2)$ between two probability distributions $\mu_1,\mu_2$ in \emph{any} nice (Polish to be exact) metric space:
\begin{equation}\label{eq:wasserstein}
  W(\mu_1,\mu_2) = \inf_{\mu\in\Gamma(\mu_1,\mu_2)}\int d(x,y)\mu(x,y)\,dx\,dy,
\end{equation} 
where $d(x,y)$ is the distance between points $x$ and $y$ in the space, and the infimum is taken over all possible transportation plans which are joint probability distributions~$\mu$ whose marginals are $\mu_1$ and $\mu_2$. If $\mu_1,\mu_2$ are represented by sand piles, then $W(\mu_1,\mu_2)$ is the minimum cost to transport pile $\mu_1$ into pile $\mu_2$, where the grain of sand at $x$ gets transported to $y$ incurring cost~$d(x,y)$.

Let $x$ and $y$ be now two points at a small distance $d_M(x,y)=\delta$ in a Riemannian manifold $M$, and let $\mu_1=\mu_x$ and $\mu_2=\mu_y$ be the normalized restrictions of the volume form in $M$ onto the balls $B_M(x,\delta),B_M(y,\delta)$ of radius $\delta$ centered at $x,y$, i.e.,
\begin{align}
  \mu_{x}(z) &= \frac{\vol(z)}{\vol[B_M(x,\delta)]}\quad\text{if } z\in B_M(x,\delta),\\
  \mu_{y}(z) &= \frac{\vol(z)}{\vol[B_M(y,\delta)]}\quad\text{if } z\in B_M(y,\delta),
\end{align}
where $\vol(z)\,dz$ is the volume element in~$M$.
Ollivier curvature between $x$ and $y$ is then defined by
\begin{equation}\label{eq:ollivier}
  \kappa_M(x,y)=1 - \frac{W_M(\mu_x,\mu_y)}{\delta}.
\end{equation}
We note that this curvature definition depends on the ball radius $\delta$ which can be any positive real number. It was shown in~\cite{ollivier2009ricci} that
\begin{equation}\label{eq:ollivier-ricci}
  \lim_{\delta\to0}\frac{\kappa_M(x,y)}{\delta^2}=\frac{\Ric(v,v)}{2(D+2)},
\end{equation}
where $D$ is $M$'s dimension, and $\Ric(v,v)$ is the Ricci curvature at $x$ along $v$ which is the unit tangent vector at $x$ pointing along the geodesic from $x$ to $y$. Ricci curvature $\Ric(v,v)$ is equal to the average of sectional curvatures at $x$ over all tangent planes containing $v$.

In words, the important convergence result in~\eqref{eq:ollivier-ricci} says that in Riemannian manifolds, the rescaled Ollivier curvature converges to Ricci curvature in the limit of small ball sizes $\delta\to0$.

It is relatively straightforward to extend the definition of Ollivier curvature from manifolds $M$ to simple unweighted graphs $G$, but there are important differences. The most crucial one is that in smooth spaces $M$ the Ollivier ball radius $\delta$ appearing in the Ollivier curvature definition~\eqref{eq:ollivier} can be any positive real number, which can tend to zero to prove the convergence in~\eqref{eq:ollivier-ricci}. In unweighted graphs $G$, however, the possible nonzero values of the distance---and so of the $\delta$---are positive integers, the hop lengths of shortest paths between pairs of nodes. The smallest possible nonzero distance is $\delta=1$, the graph distance $d_G(x,y)$ between directly connected nodes~$x$ and~$y$. Presumably because of the ``$\delta$ better be small'' logic, the Ollivier ball radius $\delta$ was routinely set to its smallest possible value $\delta=1$ in the definition of Ollivier graph curvature appearing in the past literature, even though $\delta$ in graphs can certainly be set to any positive integer value.

If $\delta=1$, then the balls $B_G(x,1)$ and $B_G(y,1)$ of radius $1$ around nodes $x$ and $y$ are the sets of the neighbors of $x$ and $y$, while the probability distributions $\mu_x$ and $\mu_y$ become the uniform distributions over these sets: $\mu_x=1/k_x$ and $\mu_y=1/k_y$, where $k_x$ and $k_y$ are the degrees of $x$ and $y$. These distributions govern the standard random walk in the graph. Thus defined, Ollivier graph curvature was shown to be bounded by $-2$ and~$1$, with the bounds achieved on infinite double-stars and complete graphs~\cite{jost2014olliviers}. The Ollivier graph curvature in this definition was measured in a great variety of synthetic and real-world networks~\cite{sandhu2015graph,ni2015ricci,sandhu2016ricci,sia2019ollivier,samal2018comparative,prokhorenkova2020global,jost2014olliviers,bhattacharya2015exact}. It was also investigated at great depths in connection to quantum gravity~\cite{trugenberger2017combinatorial,kelly2019assembly,kelly2021emergence,klitgaard2018introducing,klitgaard2018implementing,klitgaard2020quantum}.

However, it is evident that Ollivier curvature of graphs as defined above cannot converge to Ricci curvature of manifolds under any circumstances, simply because $\delta$ is no longer real but integer, so that $\lim_{\delta\to0}$ makes no sense, and because Ollivier graph curvature is always between $-2$ and $1$, while Ricci curvature of a Riemannian manifold can be any real number. Therefore, in the next section, we generalize the Ollivier graph curvature definition to make it much more versatile, capturing a much wider spectrum of graph properties that are not necessarily microscopic.

\section{Mesoscopic graph neighborhoods}\label{sec:mesoscopic}

\subsection{General idea}

The main limitation of the Ollivier graph curvature definition in the previous section is that it is ``too microscopic.'' It limits the curvature-related considerations to graph balls/neighborhoods of the radius of one hop, $\delta=1$. Such neighborhoods may be too small to ``feel'' any curvature in a general case. Recall that microscopically, any smooth curved space is locally flat. This observation instructs us to consider much larger mesoscopic graph neighborhoods with $\delta\gg1$. Unfortunately, this instruction cannot be correct either, because to talk about any curvature convergence, we really need $\delta\ll1$ since $\delta$ must go to zero in the key convergence equation~\eqref{eq:ollivier-ricci}, which is the main motivation to set $\delta$ to its smallest possible value $1$ in graphs.

The combination of the calls for $\delta\ll1$ and $\delta\gg1$ looks like a clear sign of unhealthy pathology, suggesting to abandon any hopes for Ollivier-to-Ricci curvature convergence in graphs. And indeed we believe it is obvious that there cannot be any such convergence in unweighted graphs, simply because $\delta$ cannot be less than $1$ in them.

However, these observations do not preclude the curvature convergence in weighted graphs with positive real weights. In such graphs, the Ollivier ball radius $\delta$ is back to being a positive real number, which can tend to zero because link weights can tend to zero with some characteristic rate~$\epsilon$. Suppose also that $\epsilon$ and $\delta$ tend to zero
such that $\epsilon\ll\delta$.
Then, on the one hand, we have that $\delta\to0$ as needed for convergence~\eqref{eq:ollivier-ricci}, while on the other hand, the hop-wise lengths of the radii of Ollivier's balls are $\delta/\epsilon\gg1$, so that their sizes actually grow in terms of number of nodes. This way we simultaneously satisfy the two requirements above: the Ollivier balls are small ($\delta\ll1$) and large ($\delta/\epsilon\gg1$) at the same time, hence we call them \emph{mesoscopic}. It is intuitive to think of $\epsilon$ and $\delta$ as the microscopic and mesoscopic scales of the system.

The outlined approach is very general and can be implemented in a variety of situations. In particular, the graphs do not really have to be weighted, but $\delta$ must necessarily play a role of a properly rescaled graph distance. In the next section, we make these ideas concrete in application to random geometric graphs.

\subsection{Random geometric graphs}

To describe our implementation of the general ideas from the previous subsection in application to random geometric graphs (RGGs), we first recall what RGGs are.

\begin{figure}
  \centerline{\includegraphics[width=\linewidth]{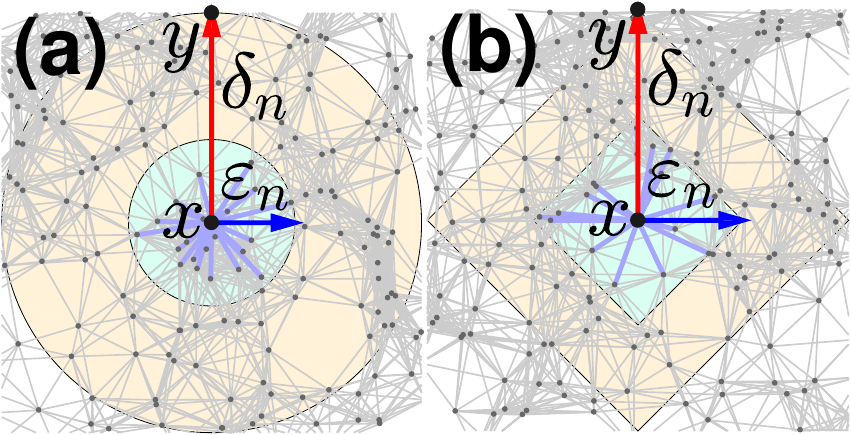}}
  \caption{{\bf Random geometric graphs in (a)~Riemannian space and (b)~Lorentzian spacetime.} The figure visualizes the flat two-dimensional Euclidian plane~(a) and Minkowski spacetime~(b), where distances are the $\ell_2$ and $\ell_1$ norms, respectively. Nodes are sprinkled randomly via the Poisson point process with rate $r=n$. The blue regions are the balls $B_M(x,\e_n)$ of radius $t=\e_n$ around node~$x$ located at the origin. In~(b) this ball is the intersection of the past light cone centered at time $t=\e_n$ with the future light cone centered at time $t=-\e_n$. Node~$x$ is connected to all other nodes in the graph that happen to lie within these blue balls. The yellow regions are the Ollivier balls $B_M(x,\delta_n)$ of radius~$\delta_n\geq\e_n$. Node~$y$ is at distance~$\delta_n$ from~$x$. Vector~$v$ (not shown) is the tangent vector at~$x$ towards~$y$. The probability distribution~$\mu_x$ is the uniform distribution over all the nodes $z$ that happen to lie within distance~$\delta_n$ from~$x$ in the graph, $z\in B_G(x,\delta_n)$. These are \emph{not} exactly the nodes lying with the yellow Ollivier balls, but at large~$n$ the difference is negligible, Appendix~\ref{ssec:aprrox-mu}.}\label{fig:fig1}
\end{figure}

Given any Riemannian manifold, which henceforth we will often call just {\it space}, the RGGs in it are defined constructively via the following two-step procedure: 1)~sprinkle points uniformly at random in the space via a Poisson point process of rate~$r>0$ driven by the volume form defined by the metric in the space~\cite{last2017lectures}, and 2)~connect by edges all pairs of points whose pairwise distances in the space are smaller than threshold $t>0$, Fig.~\ref{fig:fig1}(a). In topology, RGGs are a very fundamental object because they are $1$-skeletons of Vietoris-Rips complexes~\cite{kahle2011random} whose topology was proven to converge to the space topology under very mild assumptions~\cite{latschev2001vietoris}. In that context, our results below show that in a proper sense, the geometry of RGGs also converges to the space geometry.

We are now ready to describe our $\epsilon,\delta$-settings with RGGs in the continuum limit $n\to\infty$ suitable for curvature convergence investigations.

Given any nice $D$-dimensional Riemannian manifold~$M$, we first fix any point $x$ in it and any unit tangent vector~$v$ at $x$. We are concerned with Ricci curvature at $x$ in the $v$-direction. For a given $n$, we add the second point $y$ at distance $\delta_n$ from $x$, $d_M(x,y)=\delta_n$, along the geodesic from $x$ in the $v$-direction. The distance $\delta_n$ is also the radius of the balls in the Ollivier curvature definition above, so that it goes to zero in the limit, $\delta_n\to0$.

We then add the Poisson point process (PPP) of rate $r\sim n$ in $M$ to the two nonrandom points $x,y$. Henceforth, $f(n)\sim g(n)$ means $f(n)/g(n) \to c\in(0,\infty)$, and $f(n)\approx g(n)$ means $c=1$. The connection distance threshold $t$ in RGGs defined on top of this PPP$+\{x,y\}$ is then set to $t=\epsilon_n\to0$.
That is, all pairs of points are linked in RGG $G$ if the distance between them in $M$ does not exceed $\epsilon_n$. In the continuum limit $n\to\infty$ of this setup, our PPP samples the manifold increasingly densely (as there are on average $r\sim n$ vertices in $G$ per unit volume in $M$), while the RGGs become increasingly ``microscopic'' (as their edges span $M$'s distances no greater than $\epsilon_n\to0$).

Every created edge is then weighted as described in the next section. For any two vertices $z_1$ and $z_2$ in RGG $G$, we thus have two distances defined: the manifold distance $d_M(z_1,z_2)$ and the graph distance $d_G(z_1,z_2)$, which is the length of the shortest path between $z_1$ and $z_2$ in the \emph{weighted graph} $G$.

Ollivier's balls centered at $x,y$ are then the sets of $G$'s vertices $z$ lying within \emph{weighted graph} distances $d_G\leq\delta_n$ from vertices $x,y$: $B_G(x,\delta_n)=\{z\in G:d_G(x,z)\leq\delta_n\}$, $B_G(y,\delta_n)=\{z\in G:d_G(y,z)\leq\delta_n\}$. The probability distributions $\mu_x,\mu_y$ are the uniform distributions on the vertex sets $B_G(x,\delta_n), B_G(y,\delta_n)$:
\begin{align}
  \mu_{x}(z) &= \frac{1}{\left|B_G(x,\delta_n)\right|}\quad\text{if } z\in B_G(x,\delta_n),\label{eq:mu_G_x}\\
  \mu_{y}(z) &= \frac{1}{\left|B_G(y,\delta_n)\right|}\quad\text{if } z\in B_G(y,\delta_n).\label{eq:mu_G_y}
\end{align}
Finally, Ollivier \emph{graph} curvature $\kappa_G(x,y)$ is defined by the same equation~\eqref{eq:ollivier}, except that probability distributions $\mu_x,\mu_y$ are given by~(\ref{eq:mu_G_x},\ref{eq:mu_G_y}), Fig.~\ref{fig:fig1}(a):
\begin{equation}\label{eq:ollivier-graph}
  \kappa_G(x,y)=1 - \frac{W_G(\mu_x,\mu_y)}{\delta_n}.
\end{equation}

We reiterate that as opposed to previous literature on Ollivier graph curvature, we do \emph{not} require the RGG connection radius $\e_n$ to be equal to the Ollivier ball radius $\delta_n$, Fig.~\ref{fig:fig1}(a). However, we do not exclude this possibility either. In fact, we will see in the next section that the Ollivier-to-Ricci curvature convergence does hold even with $\e_n=\delta_n$ if links are weighted by the manifold distances between the linked nodes, and if graphs are sufficiently dense. That is, our only requirements to $\e_n$ and $\delta_n$ are that they both must go zero, and that $\e_n\leq\delta_n$.

\section{Results}\label{sec:results}

We obtain curvature convergence results for two different schemes of edge weighting in our random geometric graphs. We consider them separately in the next two sections.

\subsection{Space distance weights}

We first state our results for the settings in the previous section with the weights $w_{ij}$ of all edges $i,j$ in RGGs $G$ set equal to the distances between $i$ and $j$ in manifold~$M$:
\begin{equation}\label{eq:distance-weights}
w_{ij}=d_M(i,j).
\end{equation}

With these settings, we can show that in the continuum limit, rescaled Ollivier curvature in random \emph{graph} $G$ converges to Ricci curvature in \emph{manifold} $M$, 
\begin{equation}\label{eq:ollivier-ricci-hoorn}
    \lim_{n\to\infty} \left\langle\left|\frac{\kappa_G(x,y)}{\delta_n^2} - \frac{\Ric(v,v)}{2(D+2)}\right|\right\rangle = 0,
\end{equation}
if the RGG connectivity radius~$\e_n$ and Ollivier's ball radius~$\delta_n$ shrink to zero with~$n$ as
\begin{equation}\label{eq:delta_n}
  \e_n      \sim n^{-\alpha},\quad
  \delta_n  \sim n^{-\beta},   
\end{equation}
with exponents $\alpha,\beta$ satisfying
\begin{equation}\label{eq:weighted-conditions}
  0 < \beta \leq \alpha,\quad
  \alpha+2\beta<\frac{1}{D}.
\end{equation}
The connectivity and Ollivier balls can be equal and can shrink at the same rate,
\begin{equation}
\epsilon_n\sim\delta_n\sim n^{-\alpha},
\end{equation}
in which case \eqref{eq:weighted-conditions} becomes
\begin{equation}\label{eq:weighted-conditions-alpha}
  \alpha<\frac{1}{3D}.
\end{equation}

The expectation $\ea{\cdot}$ in~\eqref{eq:ollivier-ricci-hoorn} is w.r.t.\ the RGG ensemble, in which $\kappa_G(x,y)$ is random. The result in~\eqref{eq:ollivier-ricci-hoorn} is strong in the sense that it implies not only the convergence of the expected value of Ollivier graph curvature to Ricci curvature generalizing~\eqref{eq:ollivier-ricci} to graphs,
\begin{equation}
    \lim_{n\to\infty}\frac{\ea{\kappa_G(x,y)}}{\delta_n^2}=\frac{\Ric(v,v)}{2(D+2)},
\end{equation}
but also the concentration of random $\kappa_G(x,y)$ around its expected value,
\begin{equation}
    \lim_{n\to\infty}\mathrm{Prob}\left[\left|\frac{\kappa_G(x,y)}{\delta_n^2}-\frac{\Ric(v,v)}{2(D+2)}\right|>\e\right]=0
\end{equation}
for any $\e>0$.

We note that the conditions~(\ref{eq:weighted-conditions}) allow graphs to be arbitrarily sparse, but not exactly ultrasparse. We call graphs \emph{dense}, \emph{sparse}, and \emph{ultrasparse} (also \emph{truly sparse}) if their expected average degree
\begin{equation}
    \bar{k}_n\sim n\,\vol[B_M(\cdot,\e_n)]\sim n\e_n^D \sim n^{1-\alpha D}
\end{equation}
is $O(n)$, $o(n)$, and $\sim const.$, respectively. Observe that the exponent $\beta$ in (\ref{eq:weighted-conditions}) cannot be exactly zero because Ollivier's balls must shrink for the sake of convergence, but it can be arbitrarily close to zero. The closer the $\beta$ to zero, the closer the $\alpha$ can be to $1/D$ according to~(\ref{eq:weighted-conditions}), while $\alpha=1/D$ corresponds to the ultrasparse case with $\bar{k}_n\sim const$. We note though that the closer the $\beta$ to zero, the slower the convergence, simply because Ollivier's balls shrink too slowly in this case.

Another important observation is that the case with $\e_n=\delta_n$ corresponds to the traditional definition of Ollivier graph curvature with microscopic one-hop neighborhoods. In this case, Eq.~\eqref{eq:weighted-conditions-alpha} says that we cannot prove the convergence anywhere close to the ultrasparse limit with $\alpha=1/D$, and our simulations suggest that there is no convergence there indeed.

In general, the sparser the graphs, the more problematic the convergence becomes. This is because, the sparser the graphs, the smaller the connection radius $\e_n$. As discussed in Section~\ref{sec:mesoscopic}, if $\e_n$ is too small, and if $\delta_n=\e_n$, then all balls are of the radius of one hop, and such microscopic neighborhoods may not ``feel'' any curvature; all they can ``feel'' is locally Euclidean flatness. 

\subsubsection{Proofs}

Our proof of~\eqref{eq:ollivier-ricci-hoorn}, outlined in Appendix~\ref{sec:proof_weighted_graphs}, is not trivial. It consists of three major steps.

First, we extend the weighted shortest path distance in RGGs $G$ to space $M$. We do so by defining the following auxiliary distance in $M$. For any pair of points $x,y\in M$ we add $x$ and $y$ as nodes to the graph $G$, and connect them to all other nodes in $G$ within some distance $\lambda_n \ll \epsilon_n$ from $x,y$ in $M$. The auxiliary manifold distance between $x,y$ is then the distance between them in the extended graph. We must select $\lambda_n$ very carefully, to make sure the auxiliary distance is a sufficiently good approximation of the true distance in $M$.

Once we have this auxiliary distance in $M$, we need to approximate Ollivier's balls in RGGs $G$ with the corresponding balls in $M$. That is, on the one hand we have sets of nodes within graph distance $\delta_n$ from node $x$, while on the other hand we have sets of nodes within manifold distance $\delta_n$ from $x$. We work out this approximation by deriving an upper bound on the Wasserstein distance between the uniform probability distributions on these two balls. This allows us to work with balls defined by manifold distances, instead of graph distances.

Finally, the most nontrivial step, we must relate the discrete uniform probability distribution defined on finite sets of nodes in balls in $G$, to the corresponding nonuniform continuous distributions defined on the continuum of points in balls in $M$. To this end, using a result on matching two Poisson point processes, we first bound the Wasserstein distance between these two distributions in Euclidean space. We finally extend this result to non-Euclidean $M$, using the local flatness of the space together with the fact that the neighborhoods that we consider shrink as $\delta_n \to 0$.

Overall, at each step of this proof, we introduce an approximation error to the relation between the Wasserstein distances in graphs $G$ and manifold $M$. The crux of each step is that it must be done with care to ensure that the introduced error is upper-bounded by $\delta_n^3 \to 0$. This is because if the error is that small, then~\eqref{eq:ollivier-ricci-hoorn} follows from~\eqref{eq:ollivier-ricci}. The further details behind how we achieve this at each step of the proof can be found in Appendix~\ref{sec:proof_weighted_graphs}.

\subsubsection{Simulations}

The fact that we can prove the Ollivier$\to$Ricci convergence only under the conditions (\ref{eq:weighted-conditions}) does not mean that there is no convergence outside of this parameter region, and this is indeed what we observe in simulations.

For simulations we select the sphere, torus, and the Bolza surface as the three $2D$ manifolds of constant Ricci curvature $+1$, $0$, and $-1$. The Bolza surface~\cite{aurich1988periodic,ratcliffe2019foundations} is the simplest hyperbolic $2D$ manifold of genus $2$ with no boundaries.
We consider the Bolza surface, versus something simpler of negative curvature, such as a hyperbolic disk, because we want our manifolds to have no boundaries. We want this, because computing Ollivier curvature in simulations is a major challenge at large $n$, so that we do not want any boundary effects at any $n$.

We then fix $x$ and $y$ as described above, sprinkle $n$ points on the three manifolds uniformly at random according to the manifold volume form, and link all pairs of points at distances $\leq \e_n$ on the manifold. We then find Ollivier's balls $B_G(x,\delta_n), B_G(y,\delta_n)$, set up the uniform probability distributions $\mu_x=1/|B_G(x,\delta_n)|,\,\mu_y=1/|B_G(y,\delta_n)|$ on them, and compute all the pairwise distances $d_G(x_i,y_j)$ in graph $G$ between nodes $i\in B_G(x,\delta_n)$ and $j\in B_G(y,\delta_n)$ whose coordinates are $x_i$ and $y_j$. All this data allow us to compute the Wasserstein distance $W_G(\mu_x,\mu_y)$ by solving the linear program
\begin{equation}\label{eq:ricci_lp}
\begin{split}
W_G(\mu_x,\mu_y)&=\min_{\rho_{ij}}\sum\limits_{i,j}d_G(x_i,y_j)\rho_{ij}\mu_x,\text{ s.t.}\\
\sum\limits_i\rho_{ij}\mu_x&=\mu_y,\quad0\leq\rho_{ij}\leq 1,\quad\sum\limits_j\rho_{ij}=1,\\
\end{split}
\end{equation}
where the minimization is over transportation plans $\rho_{ij}$ whose entries describe pairwise movements of probability masses. 
Having computed this $W_G(\mu_x,\mu_y)$, Ollivier curvature is finally given by~\eqref{eq:ollivier-graph}. Further details on the simulations are in Appendix~\ref{sec:simulations}.

\begin{figure}[!t]
\centering
\includegraphics[width=1.08\linewidth]{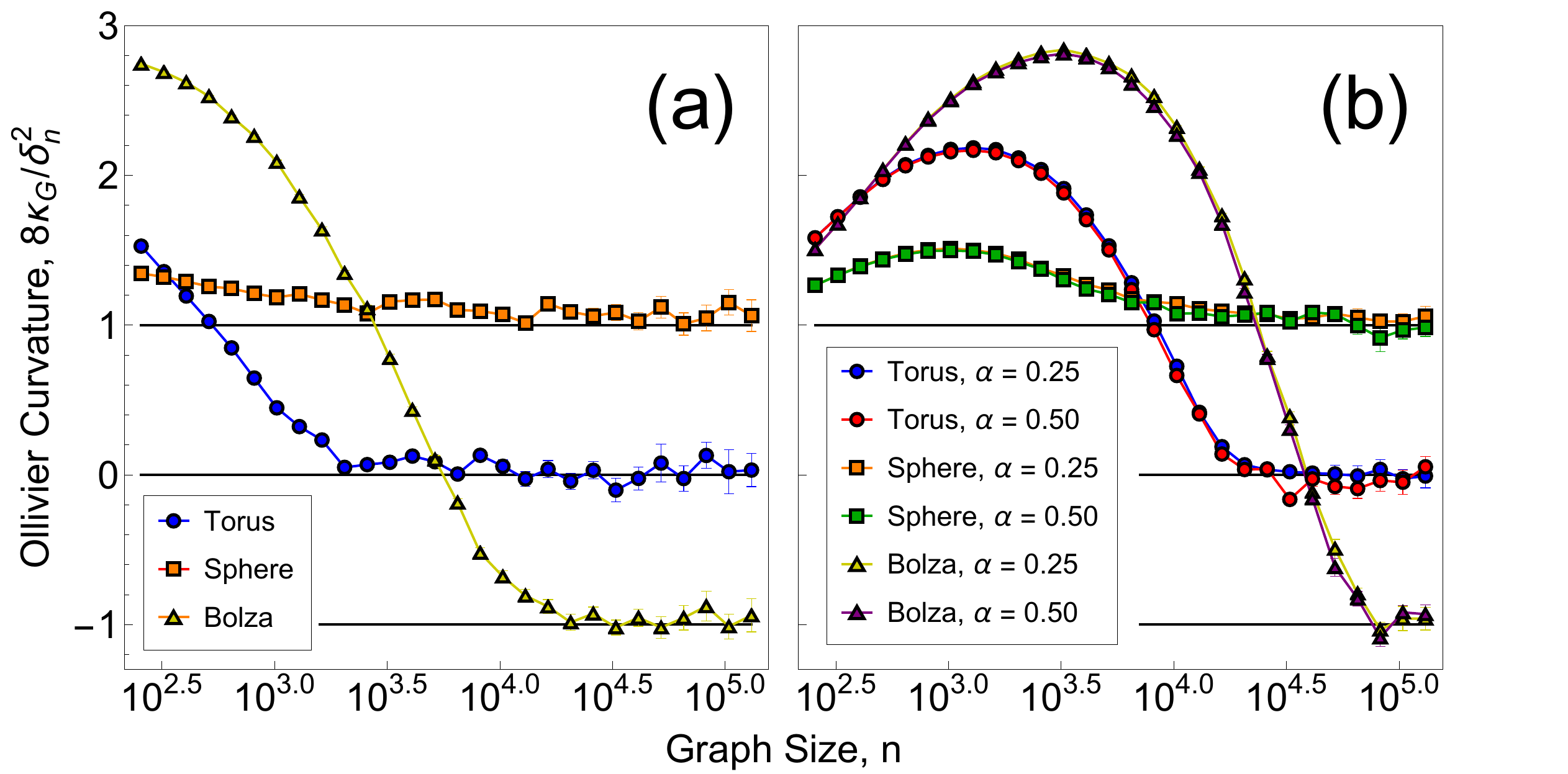}
\caption{{\bf Ollivier-Ricci curvature convergence in constant-curvature random geometric graphs.}
Panel~(a) shows the simulation data for rescaled Ollivier curvature~\eqref{eq:ollivier-graph} in \eqref{eq:distance-weights}-weighted random geometric graphs on the ($D=2$)-dimensional torus, sphere, and Bolza surface. The scaling exponents of the connectivity and Ollivier's ball radii~\eqref{eq:delta_n} in the graphs are set to $\alpha=\beta=0.16<1/6$ lying within the proof-accessible regime~(\ref{eq:weighted-conditions-alpha}).
Panel~(b) shows apparent convergence for $\alpha=\{1/4,1/2\}$ and $\beta=1/4$. Both of these two settings lie well outside of the proof-accessible regime~(\ref{eq:weighted-conditions}), while $\alpha=1/2$ corresponds to truly sparse graphs with constant average degree. The error bars represent the standard error $\sigma/\sqrt{n_s}$, where $\sigma$ is the standard deviation and $n_s=10\cdot2^{17}/n$ is the number of sampled random graphs of size $n$.
The smallest graph size $n$ is chosen to ensure the graphs are connected, while the largest size $n=2^{17}$ is bounded by memory constraints of the MOSEK package.
}
\label{fig:results}
\end{figure}

The results are shown in Fig.~\ref{fig:results}.
Remarkably, we observe that rescaled Ollivier curvature in our RGGs converges to Ricci curvature of the underlying manifold not only in the parameter regime accessible to our proofs, Fig.~\ref{fig:results}(a), but also well outside of this regime, for much sparser graphs, Fig.~\ref{fig:results}(b). In particular, we observe that the convergence is as good for truly sparse graphs with constant average degree as for denser graphs.

\subsection{$\e_n$-weights}

On the one hand, we may be not too happy that the weighted RGGs considered in the previous section contain too much information about the manifold in their edge weights $w_{ij}=d_M(i,j)$, helping the convergence ``too much.'' On the other hand, we cannot be happy at all if edges are unweighted, $w_{ij}=1$, because in this case there cannot be any curvature convergence whatsoever for the reasons discussed in Section~\ref{sec:mesoscopic}.
The maximum-happiness point appears to be the middle-ground possibility discussed in Section~\ref{sec:mesoscopic}---all edges are weighted by the same weight decaying to zero and playing the role of a system scale, {\it a la} the Planck scale.
This ``Planck scale'' in our case is the RGG connectivity radius $\e_n$.
That is, the weights $w_{ij}$ of all edges $i,j$ in our RGGs in this section are
\begin{equation}\label{eq:e-weights}
w_{ij}=\e_n.
\end{equation}

This weighting scheme is equivalent to saying that we actually work with \emph{unweighted} graphs, except that we multiply by $\e_n$ all distances in them, which are hop lengths of shortest paths. If $d_G(z_1,z_2)$ is the shortest path hop distance between vertices $z_1$ and $z_2$ in the unweighted RGG $G$, then the rescaled graphs distance $d'_G(z_1,z_2)=\e_nd_G(z_1,z_2)$ is actually a good approximation to the manifold distance $d_M(z_1,z_2)$. Unfortunately, how good this approximation is (characterized by \emph{stretch} $d'_G/d_M$), has been rigorously quantified only for the RGGs in the Euclidean plane~\cite{diaz2016relation}. Since the distances we approximate are bounded by Ollivier's ball radius $\delta_n$ that tends to zero, and since any Riemannian manifold is locally Euclidean, our results below relying on~\cite{diaz2016relation} hold for any $2D$-manifold. In fact, it turns out that the results in~\cite{diaz2016relation} are strong enough for our proof strategy in the weighted case to work essentially without a change in this unweighted case as well, Appendix~\ref{sec:proof_graphs}. While our result below holds only for $D=2$, analogous results for higher dimensions can be worked out as soon as stretch results for higher dimensional RGGs are obtained.

The result in Appendix~\ref{sec:proof_graphs} is that with weights in~\eqref{eq:e-weights}, the Ollivier$\to$Ricci convergence holds in the same strong sense~\eqref{eq:ollivier-ricci-hoorn} if
\begin{equation}\label{eq:non-weighted_conditions}
  0<\beta<\frac{1}{9},\quad
  3\beta<\alpha<\frac{1-3\beta}{2}.
\end{equation}

Similar to the weighted case, these conditions allow for arbitrarily sparse but not ultrasparse graphs ($\beta\to0$, $\alpha\to1/2$). Unlike the weighted case, these conditions imply that, as expected, the convergence is slower: $\beta<1/9$ (2D unweighted) versus $\beta<1/6$ (2D weighted). Another difference with the weighted case is that these conditions do not allow for the connectivity radius $\e_n$ to be the same as Ollivier's ball radius $\delta_n$: $\alpha\geq3\beta\Rightarrow\delta_n\gg\e_n$. We believe this is not a deficiency in our proof techniques, but a reflection of reality. It seems obvious that one must consider mesoscopic neighborhoods of size $\delta_n\gg\e_n$ in unweighted graphs since one-hop microscopic $\e_n$-neighborhoods in these graphs are too small to feel any curvature.

\section{Discussion}\label{sec:discussion}

We have established a rigorous connection between curvatures of discrete and continuous objects---random graphs and Riemannian manifolds. Anecdotally, this connection is reminiscent of what Riemann had in mind working on the foundations of Riemannian geometry~\cite{riemann1873hypotheses}.
The closest result to ours, which deals with convergence of curvature of simplicial triangulations of manifolds in Regge calculus~\cite{regge1961general,cheeger1984curvature} and its many derivatives~\cite{loll1998discrete}, is very different. Informally, while any simplex in any simplicial triangulation of any manifold is a chunk of space, graphs are much more primordial objects in that they do not have any space attached to them whatsoever.

We saw why the convergence of discrete curvature to a continuous one is such a delicate matter. If we use the standard simple definition of Ollivier graph curvature, then the immediately obvious observation is that there cannot be any convergence. We really have to consider mesoscopic regions in graphs with properly rescaled distances for the convergence to take place, an observation that has been entirely overlooked in the past. Even then the graph density and mesoscopic region sizes cannot be arbitrary, but should enter in a delicate interplay. These nontrivial delicate points are one reason why there are so few results on discrete$\to$continuous curvature convergence.

Our convergence proofs can likely be improved since our simulations suggest that the convergence holds well outside of the parameter regimes accessible by our proof techniques. Our simulations can be improved as well. We could not simulate graphs large enough to see any convergence in $\e_n$-weighted graphs, where the convergence is slow. In any case, the investigation of the exact boundaries between regimes where the convergence does and does not take place is a wide open problem, both for Ollivier curvature and for all other discrete curvature definitions mentioned in Section~\ref{sec:curvature}, for most of which it is not even known under what circumstances, if any, they do converge.

If some discrete curvature does not converge to a continuous one for some discrete structure, then it may very well happen that this discrete curvature tells us something very different about the structure from what is intended. In such a case, why should we call this ``curvature'' \emph{curvature}? This problem is not semantic but practical. If we are sure that graph curvature converges to manifold curvature, then we can be also sure that our graph curvature can be used to learn reliably manifold geometry---an exciting news for manifold learning and many other network/data science applications mentioned in the introduction. As an example of the other extreme, if we are sure that curvature of some tensor network does not converge to Ricci curvature of classical spacetime, then we can be also sure that this tensor network is not a good model of quantum spacetime.

As far as possible applications to quantum gravity are concerned, our results are an important but very high-level step that does not touch on many critical aspects. We comment on some of them.

First of all, we proved the convergence of curvature only for graph links and paths. Their curvatures first need to be averaged to result in scalar Ollivier-Ricci curvature at a vertex, and then integrated over vertex sets to yield a discrete version of the Einstein-Hilbert action. Is there any convergence in this case, and if so, then under what circumstances? Such convergence results would indicate unambiguously that the discrete Einstein-Hilbert action can be constructed from Ollivier curvature.

A much more difficult problem is the following. What we have solved is the direct problem of graph curvature convergence to space curvature when space is there. The problem in quantum gravity is actually the inverse problem of geometrogenesis~\cite{wu2015emergent,trugenberger2017combinatorial,kelly2019assembly}---there is no space, only a discrete quantum structure, which must ``look like'' space in the classical limit. That is, an illusion of continuous space must somehow emerge in the limit.

One strategy to address this problem is conceptually similar to~\cite{krioukov2016clustering}. One can consider a canonical maximum-entropy Gibbs ensemble of random graphs with fixed values of Ollivier curvature on links, paths, or vertices. The questions then are: is this ensemble equivalent (in the continuum limit) to the ensemble of random geometric graphs on a Riemannian manifold with the corresponding values of Ricci curvature, and do we have any (second-order) phase transitions in this ensemble? If the answer to the first question is \emph{yes}, then the graphs in this ensemble do look like spatial graphs, even though there is no space in their definition. If the connection between Ollivier curvature and the Einstein-Hilbert action is there, then the temperature parameter in this ensemble must be analogous to the gravitational coupling constant. In that context, the settings considered in this paper are particularly interesting. Indeed, since Ollivier curvature~\eqref{eq:ollivier-graph} is well defined for any shortest path between any pairs of nodes $x,y$ located at any distance in the graph, one can probe graph curvature at \emph{any} scale for renormalization and other purposes.

Finally, yet another critical issue is the issue of time, causality, and Lorentzian metric signature. There are two schools of thought concerning the issue. The first school believes that time and Lorentzian signature must be also emergent. This route is followed, for example, in the combinatorial quantum gravity program~\cite{trugenberger2017combinatorial,kelly2019assembly}, where time and causality are expected to appear only in the continuum limit above a certain scale, below which spacetime is Riemannian, turning into a random graph at a second-order phase transition defining quantum gravity nonperturbatively.

The second school builds in time and causality {\it ab initio}, as in causal sets or causal dynamical triangulation. For our results to be applicable in these settings, one needs first to extend them from Riemannian spaces to Lorentzian spacetimes. To the best of our knowledge, no notion of random geometric graphs or Rips complexes has been defined for Lorentzian manifolds. Since intersections of past and future light cones, Fig.~\ref{fig:fig1}(b), known as \emph{causal diamonds} or \emph{Alexandroff sets}, play the role of balls in Lorentzian geometry---specifically, the topology defined by Alexandroff sets agrees with the base topology in any nice Lorentzian spacetime~\cite{kronheimer1967structure}---we propose to define Lorentzian random geometric graphs as shown in Fig.~\ref{fig:fig1}(b). This definition is very different from causal sets, where any two time-like-separated nodes are linked, and there are no other links. In the Fig.~\ref{fig:fig1}(b) definition of Lorentzian random geometric graphs, there are both time-like and space-like links but only between nodes within causal diamonds of a finite size. Can Ollivier curvature be defined for these graphs and shown to converge to Ricci curvature of the underlying Lorentzian manifold?

\begin{acknowledgments}
We thank R.~Loll, B.~Dittrich, S.~Surya, D.~Rosset, E.~Andersen, Y.~Zhou, J.~Jost, and A.~Aranovich for useful discussions, suggestions, and comments.
This work was supported by ARO Grant Nos.~W911NF-16-1-0391 and W911NF-17-1-0491, and by NSF Grant Nos.~IIS-1741355 and DMS-1800738. Research at Perimeter Institute is supported by the Government of Canada through the Department of Innovation, Science, and Economic Development and by the Province of Ontario through the Ministry of Colleges and Universities. Experiments were conducted using the Discovery cluster, supported by Northeastern University's Research Computing team, and the B\'eluga, Graham, and Cedar clusters, supported by Compute Canada and its regional partners Calcul Qu\'ebec, Compute Ontario, and Westgrid.
\end{acknowledgments}

\onecolumngrid
\clearpage
\twocolumngrid

\appendix

\section{Ollivier$\to$Ricci convergence in distance-weighted RGGs}\label{sec:proof_weighted_graphs}

Here we outline the main steps and ingredients in our curvature convergence proof. All further lower-level details can be found in~\cite{hoorn2020ollivier}.

Let $\mu^M_x$ and $\mu^M_y$ denote the normalized restrictions of the volume form in manifold $M$ onto the balls $B_M(x,\delta_n),B_M(y,\delta_n)$ in $M$, while $\mu^G_x$ and $\mu^G_y$ denote the uniform probability distribution over nodes lying within the balls $B_G(x,\delta_n),B_G(y,\delta_n)$ in weighted RGG $G$. Observe that if
\begin{equation}\label{eq:main_part_proof_weighted}
	\langle|W_G(\mu^G_x, \mu^G_y) - W_M(\mu^M_x,\mu^M_y)|\rangle \ll \delta_n^3,
\end{equation}
then~\eqref{eq:ollivier-ricci-hoorn} follows from~\eqref{eq:ollivier-ricci}. Indeed, if \eqref{eq:main_part_proof_weighted} holds, then 
\begin{align*}
  \langle|\kappa_G(x,y)-\kappa_M(x,y)|\rangle
  &=\frac{\langle|W_G(\mu^G_x, \mu^G_y) - W_M(\mu^M_x,\mu^M_y)|\rangle}{\delta_n} \\
  &\ll\delta_n^2,
\end{align*}
so that
\[\left\langle\left|\frac{\kappa_G(x,y)}{\delta_n^2}-\frac{\kappa_M(x,y)}{\delta_n^2}\right|\right\rangle\ll1.\]

The main idea of the proof is thus to show~\eqref{eq:main_part_proof_weighted}. We proceed in three steps.

\subsection{Approximating graph distances by manifold distances}

The first difficulty we face is that the Wasserstein distances $W_G$ and $W_M$ are in different spaces with different distances: on the one hand we have RGG $G$ with weighted shortest path distance $d_G$, and on the other hand we have the manifold $M$ with distances $d_M$. How can we compare the two then? We simply extend $d_G$ to a new distance $\widetilde{d}_M$ on $M$, such that the difference between the new and original Wasserstein distances $\widetilde{W}_M$ and $W_M$ is $\ll \delta_n^3$.

This extension is accomplished by the following procedure. Let
\begin{equation}
	\lambda_n = \log(n)^{2/D} n^{-1/D},
\end{equation}
so that
\begin{equation}
	n \lambda_n^D \to \infty. \label{eq:lambda_volume_diverges}
\end{equation}
If the conditions~(\ref{eq:weighted-conditions}) hold, then we have
\begin{align}
	\lambda_n &\ll \varepsilon_n, \label{eq:lambda_bound_epsilon} \\
	\lambda_n &\ll \delta_n^3.\label{eq:lambda_bound_delta}
\end{align}

\begin{figure}
\centering
\vspace*{3pt}
\begin{tikzpicture}[scale=0.55]
	\tikzstyle{vertex}=[fill, circle, inner sep=0pt, minimum size=5pt]
	\tikzstyle{edge}=[color=black]
	
	\draw node[vertex] (u) at (0,0) {};
	\draw node[vertex] (v) at (10,0) {};
	
	\path (u)+(100:0.6) node[vertex,blue] (x0) {};
	\path (x0)+(320:2.5) node[vertex,blue] (x1) {}; 
	\path (x1)+(60:2.8) node[vertex,blue] (x2) {}; 
	\path (x2)+(340:2.9) node[vertex,blue] (x3) {};
	\path (x3)+(-30:2.5) node[vertex,blue] (x4) {};
	\path (v)+(340:0.7) node[vertex,blue] (x5) {}; 
	
	\path (u)+(-90:0.5) node {$x$};
	\path (v)+(90:0.5) node {$y$};

	\path (x0)+(90:0.7) node {\color{blue} $v_1$};
	\path (x1)+(-135:0.5) node {\color{blue} $v_2$};
	\path (x2)+(90:0.5) node {\color{blue} $v_3$};
	\path (x3)+(0:0.6) node {\color{blue} $v_4$};
	\path (x4)+(90:0.5) node {\color{blue} $v_5$};
	\path (x5)+(0:0.8) node {\color{blue} $v_6$};
	
	\draw[blue] (u) -- (x0) -- (x1) -- (x2) -- (x3) -- (x4) -- (x5) -- (v);
	
	\draw[thick] (u) circle  (1cm);
	\draw[thick] (v) circle  (1cm);
	
	\path (x4)+(-90:3) node[inner sep=0pt,minimum size=0pt] (epsilon) {};
	\path (epsilon)+(60:1.5) node {\color{blue} $\varepsilon_n$};
	\draw[dashed,blue,thick] (x4) -- (epsilon);
	
	\path (u)+(180:1) node[inner sep=0pt,minimum size=0pt] (lambda) {};
	\path (lambda)+(135:0.6) node {$\lambda_n$};
	\draw[dashed,thick] (u) -- (lambda);
	
	\draw[thick,blue] (x0) circle (3cm);
	\draw[thick,blue] (x2) circle (3cm);
	\draw[thick,blue] (x4) circle (3cm);
	
\end{tikzpicture}
\caption{Illustration of the construction of the extended graph distance $\widetilde{d}_M$. The blue points are nodes of the RGG and the two black points are the selected points on the manifold. The blue circles indicate the connection radius $\varepsilon_n$, while the black circles have radius $\lambda_n$. In this example $x$ is connected to node $v_1$ and $y$ to node $v_6$.}
\label{fig:extended_graph_distance}
\end{figure}
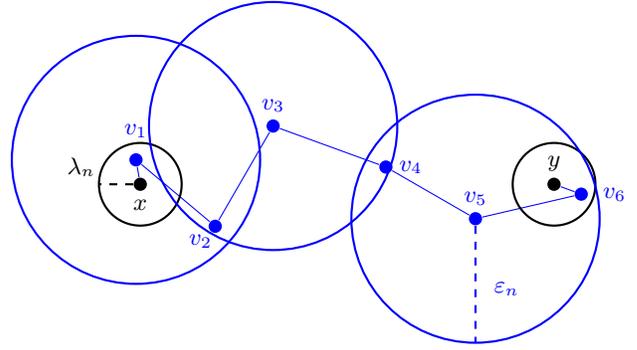

Given any two points $x,y$ in $M$ we extend the RGG $G$ by adding $x$ and $y$ to the set of $G$'s nodes, and connecting both $x$ and $y$ to all other nodes in $G$ whose distance from $x$ and $y$, respectively, in the manifold is less than $\lambda_n$, Fig.~\ref{fig:extended_graph_distance}. We then define the new manifold distance $\widetilde{d}_M(x,y)$ as the weighted shortest path distance in the extended graph. This new manifold distance $\widetilde{d}_M$ is well defined, with probability converging to one, because the expected number of points of a PPP of rate $n$ in a ball of radius $\lambda_n$ is $n\,\vol[B_M(x,\lambda_n)]\sim n \lambda_n^D$, so that the probability to find a graph vertex within distance $\lambda_n$ from any point $x$ in $M$ is given by $1-e^{-n\,\vol[B_M(x,\lambda_n)]}$ which converges to one thanks to~\eqref{eq:lambda_volume_diverges}.
Observe that if $x$ and $y$ happen to be already nodes in $G$, then the new manifold distance between them is the graph distance, $\widetilde{d}_M(x,y) = d_G(x,y)$.

\begin{figure}
\centering
\begin{tikzpicture}[scale=0.55]
	\tikzstyle{vertex}=[fill, circle, inner sep=0pt, minimum size=5pt]
	\tikzstyle{edge}=[color=black]
	
	\draw node[vertex] (u) at (0,0) {};
	\draw node (z1) at (3,0) {};
	\path (z1)+(0,0) node[vertex,red] (c2) {};
	\draw node (z4) at (9,0) {};
	\path (z4)+(0,0) node[vertex,red] (c4) {};
	\draw node[vertex] (v) at (12,0) {};
	
	\path (u)+(100:0.8) node[vertex,blue] (x0) {};
	\path (z1)+(210:0.8) node[vertex,blue] (x1) {}; 
	\path (x1)+(15:2.5) node (x31) {};
	
	\path (z4)+(300:0.9) node[vertex,blue] (x4) {};
	\path (x4)+(150:2.5) node (x32) {};
	\path (v)+(50:0.6) node[vertex,blue] (x5) {}; 
	
	\path (u)+(180:0.5) node {$x$};
	\path (v)+(0:0.5) node {$y$};
	
	\path (c2)+(90:0.5) node {\color{red} $c_{2}$};
	\path (c4)+(90:0.5) node {\color{red} $c_{m}$};

	\path (x0)+(90:0.4) node {\color{blue} $v_{1}$};
	\path (x1)+(-90:0.4) node {\color{blue} $v_{2}$};
	\path (x4)+(300:0.6) node {\color{blue} $v_{m}$};
	\path (x5)+(45:0.6) node {\color{blue} $v_{m+1}$};
	
	\draw[thick] (u) -- (z1);
	\draw node at (6,0) {$\dots\dots$};
	\draw[thick] (z4) -- (v);
	\draw[blue] (u) -- (x0) -- (x1) -- (x31);
	\draw[blue] (x32) -- (x4) -- (x5) -- (v);
	
	\draw[dashed,red] (u) circle  (1cm);
	\draw[dashed,red] (c2) circle  (1cm);
	\draw[dashed,red] (c4) circle  (1cm);
	\draw[dashed,red] (v) circle  (1cm);
	
\end{tikzpicture}
\caption{Approximating the manifold distance by the extended weighted shortest path graph distance.}
\label{fig:cover_geodesic_xy}
\end{figure}
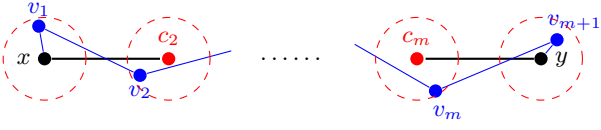

We now have to show that for any two points $x$ and $y$ in $M$, the distance $\widetilde{d}_M(x,y) = d_M(x,y) + o(\delta_n^3)$. As evident from Fig.~\ref{fig:extended_graph_distance}, $\widetilde{d}_M(x,y)\geq d_M(x,y)$. For the upper bound, let us partition the geodesic between $x$ and $y$ into $m = \lceil 3d_M(x,y)/\varepsilon_n \rceil$ sections of equal length, and let $c_1 = x, c_2, \dots, c_{m}, c_{m+1} = y$ denote the $m+1$ endpoints of this partition, Fig.~\ref{fig:cover_geodesic_xy}. Consider the $m+1$ balls $B_M(c_i, \lambda_n/4)$. The probability that every ball contains a graph vertex converges to 1. Denote those $m+1$ vertices by $v_1, \dots, v_{m+1}$.
We then observe that:
\begin{enumerate}
  \item the node $v_1$ is at distance at most $\lambda_n/2$ from $x$, and so is $v_m$ from $y$;
  \item the distance between each consecutive pair of nodes $v$ is bounded by
\[
	d_M(v_i,v_{i+1}) \le d_M(x,y)/m + \lambda_n \le \varepsilon_n/3 + \lambda_n;
\]  
  \item $d_G(v_i,v_{i+1}) \le d_M(c_i,c_{i+1}) + \lambda_n$.
\end{enumerate}
From these observations and~\eqref{eq:lambda_bound_epsilon}, we conclude that each consecutive pair of nodes $v$ is connected in the RGG and 
\begin{align*}
	\widetilde{d}_M(x,y) &\le \lambda_n + \sum_{i = 1}^{m} d_M(v_i,v_{i+1}) \\
	&\le \lambda_n + d_M(x,y)(1 + 3\lambda_n/\varepsilon_n).
\end{align*}
If $d_M(x,y) \le C \delta_n$, for some $C > 0$, then
\begin{align}
	|d_M(x,y) - \widetilde{d}_M(x,y)| &\le d_M(x,y) \frac{3 \lambda_n}{\varepsilon_n} + \lambda_n \label{eq:distance_approx} \\
	&\le 3C \lambda_n \delta_n \varepsilon_n^{-1} +\lambda_n \\
	&\sim \log(n)^{\frac{2}{D}} n^{-\frac{1}{D} - \beta + \alpha} + \lambda_n. 
\end{align}
The second term in the last equation is $\ll \delta_n^3$ by~\eqref{eq:lambda_bound_delta}, and so is the first term if~\eqref{eq:weighted-conditions} holds.
We therefore conclude that
\[
	\widetilde{d}_M(x,y) = d_M(x,y) + o(\delta_n^3)
\]

Finally, the Wasserstein distances $W_M$ and $\widetilde{W}_M$ are defined by distances $d_M$ and $\widetilde{d}_M$, respectively. Therefore we immediately conclude that
\[
	|W_M(\mu_1, \mu_2) - \widetilde{W}_M(\mu_1, \mu_2)| \ll \delta_n^3
\]
for any two probability distributions $\mu_1$ and $\mu_2$ on $M$. If now $x,y$ are also any two nodes of our RGG $G$, then $d_G(x,y) = \widetilde{d}_M(x,y)$, so that it follows that if $\mu_1$ and $\mu_2$ are defined on $G$'s nodes, then
\begin{equation}\label{eq:W_M-W_G}
	|W_M(\mu_1, \mu_2) - W_G(\mu_1, \mu_2)| \ll \delta_n^3.
\end{equation}

The calculations above thus show that from now on we can always work only with the Wasserstein distance $W_M$ in the space $M$, and this is indeed what we do at the next steps.

\subsection{Approximating probability distributions in graphs}\label{ssec:aprrox-mu}

The next hurdle we need to overcome is comparing the uniform probability distribution $\mu_x^G$ on $B_G(x,\delta_n)$ with $\mu^M_x$, the normalized restriction of the volume form to $B_M(x,\delta_n)$. The main difficulty here lies in that $d_M(x,z) \le \delta_n$ does not necessarily imply that $d_G(x,z) \le \delta_n$. In general, the shape of the intersection $B_G(x,\delta_n) \cap B_M(x,\delta_n)$ can be highly non-trivial. We tame this shape by introducing a new probability distribution $\hat{\mu}^G_x$ which is the uniform probability distribution on the set $B_G \coloneqq B_M(x,\delta_n) \cap G$ which is the set of $G$'s nodes that happen to lie within distance $\delta_n$ from $x$ in the space, i.e., all nodes in the rose ball in Fig.~\ref{fig:fig1}(a). Observe that $B_G(x,\delta_n)\subseteq B_G$ since $d_G(x,y)\geq d_M(x,y)$. The main goal of this step is to show that $\langle W_M(\mu^G_x, \hat{\mu}^G_x) \rangle \ll \delta_n^3$.

First, consider~\eqref{eq:distance_approx}. If the conditions~(\ref{eq:weighted-conditions}) hold, then $\lambda_n \ll \delta_n^3$ and $\lambda_n/\varepsilon_n \ll \delta_n^2$. It follows then that there exists a $\xi_n \ll \delta_n$ such that
\begin{equation}\label{eq:distance_approx_general}
	|d_M(x,y) - d_G(x,y)| \le d_M(x,y) \xi_n^2 + \xi_n^3.
\end{equation}
For instance, we can take the $\xi_n$ to be 
\[
	\xi_n = \max\{ \sqrt{3\lambda_n/\varepsilon_n}, \lambda_n^{1/3}\}.
\]
Define now the new radii $\delta_n^\pm = (\delta_n \pm \xi_n^3)/(1 \mp \xi_n^2)$, and note that $\delta_n^- < \delta_n < \delta_n^+$, so that
\[
	B_M(x, \delta_n^-) \subset B_M(x, \delta_n) \subset B_M(x, \delta_n^+).
\]
Let $B_G^\pm \coloneqq B_M(x, \delta_n^\pm) \cap G$ and $\mu^\pm_x$ denote the uniform probability distribution on $B_G^\pm$.
We then obtain an upper bound for $W_M(\hat{\mu}^G_x,\mu^+_x)$ by considering the following joint probability distribution a.k.a.\ transport plan $\mu$ for $(u,v) \in B_G \times B_G^+$
\[
	\mu(u,v) = \begin{cases}
		\frac{1}{|B_G^+|},			&\mbox{if } u = v\in B_G,\\
		\frac{1}{|B_G| |B_G^+|},	&\mbox{if } u\in B_G, v \in B_G^+ \setminus B_G,\\
		0	&\mbox{otherwise.}
	\end{cases}
\]
Since the Wasserstein distance $W_M(\hat{\mu}^G_x,\mu^+_x)$ is the infimum over all joint distributions~\eqref{eq:wasserstein}, and since $\hat{\mu}_x^G$ and $\mu_x^+$ are discrete, it follows that
\begin{align*}
	W_M(\hat{\mu}^G_x,\mu^+_x) &\le \sum_{u \in B_G, v \in B_G^+} d_M(u,v) \mu(u,v)\\
	&= \sum_{u \in B_G, v \in B_G^+ \setminus B_G} \frac{d_M(u,v)}{|B_G| |B_G^+|}\\
	&\le 2\delta_n^+ \frac{|B_G^+ \setminus B_G|}{|B_G^+|}.
\end{align*}
Since we are dealing with a PPP, both $|B_G^+ \setminus B_G|$ and $|B_G^+|$ are Poisson distributed with diverging means. Although they are not completely independent, it can be shown that
\[
	\left\langle \frac{|B_G^+ \hspace{-2pt} \setminus \hspace{-2pt} B_G|}{|B_G^+|} \right\rangle 
	\sim \frac{\langle |B_G^+ \hspace{-2pt} \setminus \hspace{-2pt} B_G| \rangle}
		{\langle |B_G^+ \hspace{-2pt} \setminus \hspace{-2pt} B_G| \rangle + \langle |B_G| \rangle}
	\sim \frac{(\delta_n^+)^D - (\delta_n^-)^D}{\delta_n^D}.
\]
This fraction is of the order $\xi_n^2 \ll \delta_n^2$, so that we conclude
\[
	\langle W_M(\hat{\mu}^G_x,\mu^+_x) \rangle \ll \delta_n^+ \delta_n^2 \sim \delta_n^3.
\]

Finally,
let node $z$ lie in the manifold ball of radius $\delta_n^-$ around node $x$, $d_M(x,z) \le \delta_n^-$. Then \eqref{eq:distance_approx_general} implies that the graph distance between $x$ and $z$ is upper bounded by $\delta_n$:
\[
	d_G(x,z) \le d_M(x,z)(1+\xi_n^2) + \xi_n^3 \le \frac{\delta_n - \xi_n^3}{1 + \xi_n^2} (1+\xi_n^2) + \xi_n^3 = \delta_n.
\]
This means that $z \in B_G(x,\delta_n)$, so that we have the following sandwich:
\[
	B_G^- \subseteq B_G(x, \delta_n) \subseteq B_G^+.
\]
Using calculations similar to the ones above for this sandwich, we obtain
\[
	\langle W_M(\mu^G_x,\mu^+_x) \rangle \ll \delta_n^3.
\]
Therefore, by the triangle inequality of the $W_M$ distance,
\begin{equation}\label{eq:mu^G.vs.mu^Ghat}
	\langle W_M(\mu^G_x, \hat{\mu}^G_x) \rangle
	\le \langle W_M(\mu^G_x, \mu^+_x) \rangle
	+ \langle W_M(\mu^+_x, \hat{\mu}^G_x) \rangle \ll \delta_n^3.
\end{equation}

The calculations above thus show that from now on we can work with the probability distribution $\hat{\mu}_x^G$ instead of $\mu_x^G$. The former is more convenient to work with than the latter because it is over $G$'s nodes that happen to lie within manifold ball $B_M(x, \delta_n)$, the rose ball in Fig.~\ref{fig:fig1}(a).

\subsection{Going from discrete to continuous probability distributions}\label{ssec:transition_to_manifold}

The final step is to go from the discrete probability distribution $\hat{\mu}_x^G$ to the continuous distribution $\mu_x^M$. Similarly to the goal of the previous section, the task is to show that 
\begin{equation}\label{eq:euclidean_to_manifold}
	\langle W_M(\hat{\mu}_x^G, \mu_x^M) \rangle \ll \delta_n^3.
\end{equation}
This is done by applying results on the matching distance between a PPP and points of a grid on the same space~\cite{leighton1989tight,shor1991minimax,talagrand1992matching,talagrand1994matching}. 

Note that both $\hat{\mu}_x^G$ and $\mu_x^M$ are now defined on $B_M(x,\delta_n)$. We now need to devise a transport plan from $\hat{\mu}_x^G$ to $\mu_x^M$ that assigns to any (measurable) set $A \subseteq B_M(x,\delta_n)$ how much mass from each PPP point in $B_G=G\cap B_M(x,\delta_n)$ is used to make up $\mu_x^M(A)$. We do so as follows.

First, suppose that space $M$ is flat, the $D$-dimensional Euclidean space. We construct our transport plan in three steps:
\begin{enumerate}
  \item we first place a grid on the space;
  \item we next find a minimal matching between this grid and the PPP points in $B_G$; and
  \item finally, for each (measurable) set $A$ in $B_M(x,\delta_n)$, we find all the PPP points in $B_G$ that are matched to the grid points that lie in $A$.
\end{enumerate}
Observe that this transport plan is such that every PPP point, i.e., graph vertex in $B_G$ contributes an equal fraction of its mass to make up $\mu_x^M(A)$.
This plan is also such that the largest distance any amount of mass has to move is given by the largest matching distance. Relying on the results from~\cite{leighton1989tight,shor1991minimax,talagrand1992matching,talagrand1994matching}, it can then be shown that
\begin{equation}
	\langle W_M(\hat{\mu}_x^G, \mu_x^M) \rangle \sim \log(n) n^{-1/D} \ll \delta_n^3,
\end{equation}
where the last inequality holds if the conditions~(\ref{eq:weighted-conditions}) are satisfied.

Next, we relax the condition that the space $M$ is flat. Indeed, it can be any nice Riemannian manifolds because $\delta_n \to 0$, so that we can map the ball $B_M(x,\delta_n)$ to the flat $D$-dimensional tangent space at $x$ using the exponential map. We have to be careful here though, because the exponential map does not preserve distances. Still, by fixing a sufficiently small neighborhood $U$ around the origin of the tangent space, we can ensure that for large enough $n$, $\exp^{-1} B_M(x,\delta_n) \subset U$, and the distances are distorted by a fixed small amount:
\[
	B_D\left(0, \frac{\delta_n}{1 + \xi}\right) \subseteq \exp^{-1} B_M\left(x,\delta_n\right)
	\subseteq B_D\left(0, \frac{\delta_n}{1 - \xi}\right) \subset U,
\] 
where $B_D(0,\delta)$ is the $D$-dimensional Euclidean ball of small radius $\delta$ around the origin in the tangent space.

Finally, relying on this mapping and on the PPP Mapping Theorem~\cite{last2017lectures} that applies to~$U$, we see that the mapped PPP is still a PPP with intensity $\sim n$, now relative to the Euclidean volume form. Moreover, since the mapped ball is sandwiched between two balls whose radii scale as $\delta_n$, the results for the $D$-dimensional Euclidean space above yield matching lower and upper bounds for the Wasserstein distance, from which~\eqref{eq:euclidean_to_manifold} follows.

The combination of \eqref{eq:euclidean_to_manifold}, \eqref{eq:mu^G.vs.mu^Ghat}, and \eqref{eq:W_M-W_G} yields~\eqref{eq:main_part_proof_weighted}, completing the proof.

\section{Ollivier$\to$Ricci convergence in $\e_n$-weighted RGGs}\label{sec:proof_graphs}

We first observe that the last step in the proof outlined in Appendix~\ref{sec:proof_weighted_graphs} does not rely on any graph distances at all. Further, the key element in the second step is~\eqref{eq:distance_approx_general},
which is also the fundamental ingredient in the first step in showing that distance $\widetilde{d}_M$ is a good approximation to distance $d_M$, leading to~\eqref{eq:W_M-W_G}.
Therefore, the convergence proof for $\e_n$-weighted graphs is complete as soon as \eqref{eq:distance_approx_general} is established for the new graph distance $d_G^\prime = \varepsilon_n d_G$, where $d_G$ is now the vanilla shortest path (hop count) distance in graph $G$. We also note that we have to do this only for the $2$-dimensional Euclidean space, since the neighborhoods shrink and we can apply techniques similar to the ones in Appendix~\ref{ssec:transition_to_manifold} to extend the result to any curved manifold.

We first observe that the conditions~(\ref{eq:non-weighted_conditions}) imply that $\varepsilon_n \gg \sqrt{\log(n)/n}$. In this case, the main stretch result from~\cite{diaz2016relation} reads:
\begin{equation}\label{eq:dG'-dM}
	|d_G^\prime(x,y) - d_M(x,y)| \le d_M(x,y) \gamma_n + \varepsilon_n,
\end{equation}
where $\gamma_n = \max\{\gamma_n^{(1)}, \gamma_n^{(2)}, \gamma_n^{(3)}\}$, and where the three $\gamma_n$s have somewhat complicated explicit expressions with different constants and different scalings with $n$ which can be shown to be
\begin{align*}
	\gamma_n^{(1)} \hspace{-2pt} \sim \hspace{-2pt} 
		\left(\frac{\log{n}}{n\varepsilon_n^2}\right)^{\frac{2}{3}}\hspace{-4pt}, \quad
	\gamma_n^{(2)} \hspace{-2pt} \sim \hspace{-2pt} 
		\left(\frac{\log{n}}{n \varepsilon_n^2}\right)^2\hspace{-4pt}, \text{ and }
	\gamma_n^{(3)} \hspace{-2pt} \sim \hspace{-2pt} 
		\left(\frac{1}{n \varepsilon_n^2}\right)^{\frac{2}{3}}\hspace{-4pt},
\end{align*}
if the conditions~(\ref{eq:non-weighted_conditions}) hold. Since these conditions also imply that $n \e_n^2 \to \infty$, it follows that
\begin{equation}
	\gamma_n \sim \left(\frac{\log{n}}{n\varepsilon_n^2}\right)^{\frac{2}{3}}
\end{equation}
in our case.

The $\e_n$ in \eqref{eq:dG'-dM} is $\ll \delta_n^3$ if the lower bound in the second condition in~\eqref{eq:non-weighted_conditions} holds. Therefore, to translate~\eqref{eq:dG'-dM} to~\eqref{eq:distance_approx_general}, all we need to do is to show that $\gamma_n \ll \delta_n^2$. If this holds, then we can simply take the $\xi_n$ in \eqref{eq:distance_approx_general} to be $\xi_n = \max\{\sqrt{\gamma_n}, \varepsilon_n^{1/3}\}$.
The requirement $\gamma_n \ll \delta_n^2$ is satisfied if
\begin{equation}
  \frac{\log{n}}{n\e_n^2\delta_n^3} \sim \log(n)\,n^{2\alpha+3\beta-1}\ll 1,
\end{equation}
which is indeed true if the upper bound in the second condition in~\eqref{eq:non-weighted_conditions} holds.

\section{Numerical methods}\label{sec:simulations}
The computation of Ollivier graph curvature in simulations consists of three tasks: graph construction, distance matrix computation, and Wasserstein distance computation. We rely on a hybrid multi-core solution which parallelizes each of these tasks, with the first two utilizing custom CUDA methods split among multiple GPUs, and the last utilizing multiple CPU cores via the MOSEK software package.

\subsection{Graph construction}

We construct random geometric graphs of size $n$ by first sampling $n$ points uniformly at random, according to the manifold volume form, on the surface of a unit torus, sphere, and Bolza surface. This is trivial for the first two manifolds, and less so for the Bolza surface.

The simplest representation of the Bolza surface~\cite{aurich1988periodic,ratcliffe2019foundations} is the hyperbolic octagon with vertices at complex coordinates
\begin{equation}
  o_k=R\exp\left[\frac{i\pi}{4}\left(\frac{1}{2}+k\right)\right],\quad k=0,1,\ldots,7,
\end{equation}
where $R=2^{-1/4}$ in the Poincar\'e disk model of the hyperbolic plane. When the opposite sides of this octagon are glued, the Bolza surface of constant curvature $-1$ is formed. To sprinkle points uniformly at random onto this octagon, we first sprinkle them uniformly at random, according to the volume form in the Poincar\'e model, which in the polar coordinates is
\begin{equation}
  dV =\frac{4r\,dr\,d\theta}{(1-r^2)^2},
\end{equation}
onto the Poincar\'e disk of radius $R$, and then remove those points that lie in this disk but do not lie in the octagon. We do this removal in the Klein disk model of the hyperbolic plane, because geodesics---octagon sides, in particular---are straight lines there. The map to go from the Poincar\'e $(r,\theta)$ to Klein $(r_K,\theta)$ polar coordinates is $(r_K,\theta)=(2r/(1+r^2),\theta)$, and in the latter coordinates, the coordinates of points that lie in the octagon meet the condition $r_K<r_c[\phi(\theta)]$, where
\begin{align}
r_c(\phi)&=R_K\frac{\cos(\pi/8)}{\cos(\pi/8-\phi)},\quad\phi\in[0,\pi/4],\\
\phi(\theta)&=\left(\theta-\frac{\pi}{8}(1+2k(\theta))\right)\mod 2\pi,\\
k(\theta)&=\left\lfloor\frac{4}{\pi}\left(\theta-\frac{\pi}{8}\right)\right\rfloor\mod 8,
\end{align}
where $R_K=2R/(1+R^2)=2^{5/4}/(1+\sqrt{2})$ is the octagon radius in the Klein model. The nodes that do not pass this test are thrown out. The remaining nodes lie in the octagon, and we work back with their Poincar\'e coordinates, $(r,\theta)=(1/r_K-\sqrt{1/r_K^2-1},\theta)$.

After all the $n$ nodes are sprinkled, we add two additional nodes $x,y$ separated by distance $\delta_n$, which are the centers of the two balls $B_G(x,\delta_n),B_G(y,\delta_n)$.
Given the coordinates of all the $n+2$ nodes, all node pairs are linked whenever their pairwise distance on the surface is below the connection threshold $\epsilon_n$. The distances on the torus and the sphere are 
\begin{align}
d_T(x_1,y_1;x_2,y_2)&=\left[\left(\frac{1}{2}-\left\lvert\frac{1}{2}-\left\lvert x_1-x_2\right\rvert\right\rvert\right)^2\right.\\
&+\left.\left(\frac{1}{2}-\left\lvert\frac{1}{2}-\left\lvert y_1-y_2\right\rvert\right\rvert\right)^2\right]^{1/2},\nonumber\\
d_S(\theta_1,\phi_1;\theta_2,\phi_2)&=\arccos\left[\cos\theta_1\cos\theta_2\right.\\
&\left.+\sin\theta_1\sin\theta_2\cos(\phi_1-\phi_2)\right].\nonumber
\end{align}
On the Bolza surface, the distance formula is more complicated. We describe it next.

The Bolza surface can also be considered as the factor space $\HH^2/\mathcal{F}$~\cite{ratcliffe2019foundations}, where $\HH^2$ is the hyperbolic plane, and the Fuchsian group $\mathcal{F}$ is defined by its eight generators
\begin{equation}
g_k=\begin{pmatrix}a&b e^{ik\pi/4}\\
be^{-ik\pi/4}&a
\end{pmatrix},\quad k=0,1,\ldots,7,
\end{equation}
where $a=1+\sqrt{2}$ and $b=\sqrt{a^2-1}$.
This group acts on the hyperbolic plane by linear fractional transformations.
Let $\mathcal{N}$ be the set of the 48 elements of $\mathcal{F}$ given by
\begin{equation}
\hat g_{k,\ell}=\prod\limits_{t=0}^\ell g_{k+3t},
\end{equation}
where $k=0,1,\ldots,7$ and $\ell=0,1,\ldots,5$. One can check that these elements map the original octagon to its 48 either side- or vertex-adjacent octagons in the tessellation of the hyperbolic plane induced by these identical octagons. Define $\mathcal{N}_0\coloneqq\mathcal{N}\cup I$ to be these 48 elements plus the identity element $I$. One can show that the distance on the Bolza surface between two points with complex coordinates $z_1,z_2$ ($z=re^{i\theta}$) is given by
\begin{align}
  d_B(z_1,z_2)&=\min_{\hat g\in\mathcal{N}_0} d_H(z_1,\hat g z_2),\text{ where}\\
  d_H(z_1,z_2)&=2\arctanh\left|\frac{z_1-z_2}{1-z_1^*z_2}\right|,
\end{align}
is the hyperbolic distance between the points in the Poncar\'e model.

We compute what nodes are linked in parallel using a multi-GPU solution, where each thread works on one pair of nodes, except in the case of the Bolza surface where each thread computes one or two of 49 distances and each thread block works on two node pairs. The adjacency matrix is first tiled such that (1) the data required to generate each tile fits on a single GPU and (2) the total number of tiles is a multiple of the number of GPUs available. This decomposition provides a scalable solution independent of the number or type of GPUs present in the system. The algorithm is improved by using the shared L1 memory cache to store the Bolza generators and partial results, warp shuffling to accumulate results, function templating to eliminate kernel branches, and asynchronous CUDA calls to pipeline data transfers across different GPUs.

\subsection{Distance matrix computation}

After constructing the graph, we identify the balls $B_G(x,\delta_n)$ and $B_G(y,\delta_n)$. When $\epsilon_n=\delta_n$, these are simply the nearest neighbors of $x$ and $y$, respectively; otherwise, they are calculated using Dijkstra's algorithm~\cite{dijkstra1959note}. For the largest graphs we consider, the size of these balls can be up to tens of thousands of nodes, resulting in over billions of pairwise distances between the nodes in the two balls. The overall simulation bottleneck is the computation of the matrix of the weighted shortest-path distances between these large sets of nodes.

We employ a custom multi-GPU $A^*$ search algorithm using the methods described in~\cite{zhou2015massively}. The standard $A^*$ algorithm works by constructing a priority queue of visited nodes $z$ using a binary heap. In computing the shortest path distance between a source node $x_i$ and destination node $y_j$, priorities are assigned to node $z$ in the queue according to the heuristic function
\begin{equation}
f(z)=d_G(x_i,z)+h(z,y_j),
\end{equation}
where $d_G(x_i,z)$ is the weighted graph distance between $z$ and source $x_i\in B_G(x,\delta_n)$, while $h(z,y_j)$ is the lower-bound estimate of the distance between $z$ and destination $y_j\in B_G(y,\delta_n)$.
As soon as new node $z'\in B_G(y,\delta_n)$ is added to the priority queue, its weighted graph distance $d_G(x_i,z')$ is added to the distance matrix. 

The $A^*$ algorithm is implemented on a single GPU by utilizing multiple priority queues, one per CUDA thread, so that it is efficient for graphs with large average degrees of the order of hundreds. Each priority queue extracts multiple states, after which we detect duplicates using a technique called parallel hashing with replacement, which is a modification of the cuckoo hashing scheme that avoids hash conflicts by allowing some duplicates to remain~\cite{zhou2015massively,pagh2001cuckoo}. During this step, the heuristic for extracted nodes is updated, and they are re-added to the priority queues using a parallel heap insertion. This procedure continues until the destination node $y_j$ has been extracted by at least one of the priority queues.

\subsection{Wasserstein distance computation}

Having in place the distance matrix between the sets of nodes in the two balls $B_G(x,\delta_n)$ and $B_G(y,\delta_n)$, we compute the Wasserstein distance between the uniform probability distributions on these balls by solving the linear program~\eqref{eq:ricci_lp}. For the simulations presented here, we found it sufficient to use the MOSEK package~\cite{andersen2000mosek} as long as the number of variables, given by $|B_G(x,\delta)||B_G(y,\delta)|$, is roughly less than $1.3\times10^9$, past which we run out of memory. In most cases with smaller balls, however, we solve the linear program quite quickly using the standard primal-dual interior point method~\cite{mehrotra1992implementation}.

\onecolumngrid
\clearpage
\twocolumngrid

\onecolumngrid

\end{document}